\documentclass[twoside,12pt]{article}

\usepackage{amsmath}
\usepackage{amssymb}
\usepackage{amscd}
\usepackage{amsthm}

\numberwithin{equation}{section}

\theoremstyle{plain}
\newtheorem{theorem}{Theorem}[section]

\theoremstyle{remark}

\theoremstyle{definition}

\newcommand{\R}{\mathbb R}
\newcommand{\C}{\mathbb C}
\newcommand{\h}{\mathbb H}

\def\ra{\rightarrow}

\def\e{\emph}

\def\b{\begin}

\begin{document}

\title{    
{\bf{Some examples of quasiisometries of nilpotent Lie  groups}} }
\author{Xiangdong Xie\footnote{Partially supported by NSF grant DMS--1265735.}}
\date{  }
\maketitle

\begin{abstract}
  We construct   quasiisometries   of   nilpotent Lie groups.
    In particular, for any simply connected nilpotent Lie group
       $N$, we construct quasiisometries from $N$ to itself that is
       not at finite distance from any map that is a  composition
       of left translations and automorphisms.  We also construct biLipschitz maps of the Heisenberg groups that send vertical lines to 
curves  that are   not  vertical   lines.

\end{abstract}








 {\small {\bf{Mathematics Subject
Classification (2010).}}     20F65,  22E25,  53C17,  53C23.










\setcounter{section}{0} \setcounter{subsection}{0}

\section{Introduction}\label{s0}

We construct   quasiisometries   of   simply connected   nilpotent
Lie groups.
  These examples   show that quasiisometries between nilpotent Lie groups are quite flexible.

  By Pansu's theorem \cite{P},  every quasiconformal map of  the quarternion Heisenberg group (equipped with a  Carnot metric)   must be a composition of left translations
      and automorphisms.
      In  conversation  with the author,    Yves  Cornulier    asked whether  there is 
    any simply connected nilpotent  Lie group $N$ such that 
    every   self quasiisometry of   $N$ is  at finite distance from a  map   that    is
       the composition of left translations and automorphisms. On a separate occasion, 
  Bruce Kleiner asked    if  
   the quarternion   Heisenberg group  
  is such a group.   
       In this note,   we  provide a negative answer to these questions.    In fact,
  for any simply connected nilpotent Lie group
       $N$ (equipped with a  left invariant Riemannian metric, and in the case of a Carnot group $N$  equipped  either with a Carnot metric or a left invariant Riemannian metric), we construct quasiisometries from $N$ to itself that is
       not at finite distance from any map that is a  composition
       of left translations and automorphisms.    Next we explain the construction  in the case of Carnot groups.

  Let   $N$ be a Carnot group with Lie algebra gradation  $\mathcal{N}=V_1\oplus  \cdots \oplus V_r$.  Fix any Euclidean metrics on  $V_1$ and
   $V_r$,  and     equip $N$ with the Carnot metric associated with the Euclidean metric on    $V_1$.  
       Let $g: V_1\ra V_r$ be any   large scale Lipschitz map;  
      that is, there exist  constants $L\ge 0$, $A\ge 0$ such that
    $$|g(x)-g(y)|\le L\cdot |x-y|+A$$
        for all $x, y\in V_1$.    Here $|\cdot|$ denotes the Euclidean norm.
      Let $\pi_1:   \mathcal{N}\ra V_1$ be the projection with respect to the gradation above.
         Since the exponential map of $N$ is a diffeomorphism,  we
         shall identify $N$ and its Lie algebra  $\mathcal{N}$.
             We define the shear map
             $S_g: \mathcal{N}  \ra \mathcal{N}$ by
               $$S_g(n)=n+  g(\pi_1(n)).$$
              Then $S_g$ is a quasiisometry  from $N$ to itself.  
                If  $g$ is   not at finite distance from any linear map, then $S_g$ is not at finite distance from any map
                   that is a composition of left translations and
                   automorphisms.  See     Section \ref{s1}  for more
                   details.
    We remark that  in general $S_g$  is not biLipschitz even if   $g$ is Lipschitz.

 Similar constructions can be carried out on any simply connected nilpotent Lie group equipped with left invariant Riemannian metrics.
   In this case    the  map  $S_g$ is indeed biLipschitz (with respect to the left invariant Riemannian metric)  if $g$ is Lipschitz.    We remark that  one can also consider analogous constructions 
 in the setting of finitely generated nilpotent groups.

Our construction in the case of quarternion Heisenberg groups (equipped with a Carnot metric)  yields
  quasiisometric homeomorphisms that are not at finite distance from
  any biLipschitz map.

In addition to the shear maps,   we record two other examples of
quasiisometries between nilpotent Lie groups.
 They are actually biLipschitz maps of the 3-dimensional
     Heisenberg group (equipped with a Carnot metric).

Quasiisometries between  solvable Lie groups often
    preserve certain foliations,  see \cite{EFW1}, \cite{EFW2},
        \cite{FM},
         \cite{SX}, \cite{X1}, \cite{X2}, \cite{X3}.
When studying quasiisometries between nilpotent Lie groups,
  one naturally wonders whether quasiisometries coarsely  preserve cosets of
  certain special subgroups, for example,  cosets of 
  the center.      Our  first example is a biLipschitz map of the Heisenberg group
    (equipped with a Carnot metric)   which does not coarsely preserve the foliation by vertical lines
(which are cosets of the center). 
 It maps each vertical line to within a bounded distance of a vertical line, but there is no uniform bound
  on the distance from vertical lines, 
   see   Lemma \ref{property.of.F}.    In particular, it maps vertical lines to curves that are not vertical 
lines.

The second example is a biLipschitz map of the Heisenberg group 
  (equipped with a Carnot metric) 
with the following properties:   (1) it
 maps vertical lines to vertical lines; (2) the induced map on the asymptotic cones
     is  the identity map; (3) it is
  not  at finite distance from any   map that is a composition of left
translations, automorphisms and  shear  maps.

The moral of these examples is that although  quasiisometries
between nonabelian nilpotent Lie groups  are more rigid than those
of the Euclidean spaces, they are still quite flexible.

\noindent {\bf{Acknowledgment}}. {The examples were discovered while
I was visiting MSRI in the fall of 2011. I would like to express
 my gratitude   to MSRI for its hospitality. I also would like to thank Bruce  Kleiner for
  useful discussions.}

\section{Shear maps}\label{s1}

  In this section we will define the shear maps,    and    determine
    when  they are quasiisometries  and when they   are  at finite distance  from an affine map.
        We will first consider the easier case of Carnot groups, then
    treat   the general case of simply connected
   nilpotent Lie groups.

    Since the exponential map of a
 simply connected
   nilpotent Lie group  is a diffeomorphism,   we will identify a  simply  connected
     nilpotent Lie group $N$ with its Lie algebra
      $\mathcal{N}$ via the exponential map.
     The group operation  on $\mathcal{N}$   is then given by the
       Campbell-Baker-Hausdorff (CBH) formula.
       The CBH formula expresses   the product
        $x *y$ ($x, y\in \mathcal{N}$) in  terms of  Lie brackets  involving $x$ and $y$.   The first a few terms
        are shown below:
        \begin{equation}\label{CBH}
          x*y=x+y+\frac{1}{2} [x,y]+\frac{1}{12}[x,[x,y]]-\frac{1}{12}[y, [x,y]]+\cdots.
          \tag{CBH formula}
          \end{equation}
       Notice that if  $N$ is $r$-step, then  all the $i$-fold Lie brackets  with $i\ge r$ vanish.

Let $\text{Aff}(N)$ be the group of homeomorphisms of $N$ generated
by
 left translations and  automorphisms. An element of $\text{Aff}(N)$ will
 be called  an {\it{affine map}} of $N$.  It contains the subgroup
 consisting of left translations, which can be identified with $N$
 and shall also be  denoted by $N$.
   Another subgroup is the group  $\text{Aut}(N)$ of automorphisms  of $N$.
     It is easy to check that   $N$   is normal in $\text{Aff}(N)$
     and we have  $\text{Aff}(N)=N \rtimes \text{Aut}(N) $.

  As indicated above, we first consider the Carnot group case.   So
  let   $N$ be a $r$-step  Carnot group with Lie algebra
$\mathcal{N}$.   Recall that $\mathcal N$
 can be written as a direct sum of vector subspaces
    $\mathcal
 N=V_1\oplus V_2\oplus\cdots \oplus V_r$  such that $[V_1,
V_i]=V_{i+1}$ for all $1\le i\le r$,
    where we set $V_{r+1}=\{0\}$.
    Let $\pi_i:\mathcal{N}\ra V_i$ be the projection with respect to
 the  direct sum decomposition above.
      We fix an inner product on each $V_i$.
  The   left invariant  homogeneous metric  $d_h$  is  determined by the homogeneous
  norm
    $$|n|_h=\sum_{i=1}^r |\pi_i(n)|^{\frac{1}{i}}  \;\; \text{for}\;\; n\in\mathcal{N}.$$
  So    $d_h(n_1, n_2)=|(-n_1)*n_2|_h$  for $n_1, n_2\in  \mathcal N$.  
  We remark that  $d_h$ is not a  metric in general since it  does not satisfy the triangle inequality.  It satisfies a generalized triangle inequality: there is some $M\ge 1$ such that    
 $$d_h(n_1, n_3)\le  M\cdot [d_h(n_1, n_2)+d_h(n_2, n_3)]$$
    for all $n_1, n_2, n_3\in \mathcal N$.    Let $d$  be the Carnot metric 
    on $N=\mathcal{N}$  determined   by 
    the norm on $V_1$. 
       Recall that  the  Carnot metric is biLipschitz equivalent with the
     homogeneous metric  $d_h$.

   Let $g: V_1\ra V_r$ be any  map.  
 The shear map  $S_g: (\mathcal{N}, d)\ra (\mathcal{N}, d)$     is defined    by
              $$S_g(n)=n+  g(\pi_1(n)).  $$
  A  map $g: X\ra Y$ between two metric spaces is called a 
 {\it{large scale Lipschitz map}}
  if  there exist  constants $L\ge 0$, $A\ge 0$ such that
    $$d(g(x_1),  g(x_2))\le L\cdot d(x_1,   x_2)+A\;\;\;  \text{for all}\;\;\;  x_1,   x_2\in  X.$$

\b{theorem}\label{main1}
  Let $S_g$ be the shear  map   defined above.
           Then\newline
           (1) $S_g$ is a quasiisometry  from $(\mathcal{N}, d)$ to itself  if and only if  
  $g: V_1\ra V_r$ is a large scale Lipschitz map  with respect to the Euclidean metrics on   $V_1$  and $V_r$;\newline
      (2)  Let $g$  be a large scale Lipschitz map.   Then  $S_g$ is  at finite distance from an affine  map if and
      only if    $g$ is   at finite distance from a linear map.

              \end{theorem}
\b{proof} (1)  We first assume  $g: V_1\ra V_r$ is a large scale Lipschitz map
  and shall prove  that  $S_g$ is a quasiisometry   of   $(\mathcal{N}, d)$.
  Since the Carnot metric  $d$ is biLipschitz equivalent with $d_h$,
   it suffices to prove the   claim  for   $d_h$.
  Let $n_1, n_2\in \mathcal{N}$.  Since  both
  $g(\pi_1(n_1))$  and
$g(\pi_1(n_2))$  lie in $V_r$, we have
  $$d_h(g(\pi_1(n_2)), g(\pi_1(n_1)))=|g(\pi_1(n_2))-
  g(\pi_1(n_1))|^{\frac{1}{r}}.$$
       The CBH formula yields
       $$\pi_1((-n_1)*n_2)=\pi_1(n_2)-\pi_1(n_1).$$
           By the definition
       of $d_h$ we have
         $$d_h(n_1, n_2)=d_h(0, (-n_1)*n_2)\ge |\pi_1(n_2)-\pi_1(n_1)|.$$
      As  $g$ is a $(L, A)$ large scale
  Lipschitz map, we   obtain
 \begin{align*}
d_h(g(\pi_1(n_2)), g(\pi_1(n_1)))&=|g(\pi_1(n_2))-
  g(\pi_1(n_1))|^{\frac{1}{r}}\\
 & \le
   ( L\cdot |\pi_1(n_2)-  \pi_1(n_1)|+A)^{\frac{1}{r}} \\
    &\le
   L^{\frac{1}{r}}\cdot  |\pi_1(n_2)-  \pi_1(n_1)|^{\frac{1}{r}}+
   A^{\frac{1}{r}}\\
  &   \le  L^{\frac{1}{r}}\cdot  d_h(n_1, n_2)^{\frac{1}{r}}+
   A^{\frac{1}{r}}.
\end{align*}
    Since both  $g(\pi_1(n_1))$  and  $g(\pi_1(n_2))$  lie in the
    center of
      $\mathcal{N}$,  the  CBH formula implies
  $$\{-n_1-g(\pi_1(n_1))\}*\{n_2+g(\pi_1(n_2))\}=(-n_1)*n_2*\{g(\pi_1(n_2))-g(\pi_1(n_1))\}.$$
    Hence
 \begin{align*}
  d_h(S_g(n_1), S_g(n_2))&=d_h(n_1+g(\pi_1(n_1)),  n_2+g(\pi_1(n_2)))\\
  &=d_h(0, \{-n_1-g(\pi_1(n_1))\}*\{n_2+g(\pi_1(n_2))\})\\
  &=d_h(0,  (-n_1)*n_2*\{g(\pi_1(n_2))-g(\pi_1(n_1))\})\\
  &  =d_h(g(\pi_1(n_1))-g(\pi_1(n_2))),(-n_1)*n_2)\\
  &  \le M\cdot\{d_h(g(\pi_1(n_1))-g(\pi_1(n_2)),0)+d_h(0,(-n_1)*n_2)\}\\
  &  =M\cdot\{d_h(g(\pi_1(n_2)), g(\pi_1(n_1)))+  d_h(n_1, n_2)\}\\
  &\le M\cdot\{ L^{\frac{1}{r}}\cdot  d_h(n_1, n_2)^{\frac{1}{r}}+
A^{\frac{1}{r}}   +  d_h(n_1, n_2)\}.
\end{align*}
  If  $d_h(n_1, n_2)\le 1$, then
$$d_h(S_g(n_1), S_g(n_2))\le
  M\cdot\{L^{\frac{1}{r}}+A^{\frac{1}{r}}   + 1\}.$$
    If   $d_h(n_1, n_2)\ge 1$, then
$$d_h(S_g(n_1), S_g(n_2))\le  M  (L^{\frac{1}{r}}+1) \cdot d_h(n_1,
n_2)+ M A^{\frac{1}{r}}.$$
  Hence $S_g$ is $(L', A')$ large scale Lipschitz with
  $L'=M (L^{\frac{1}{r}}+1)$ and $A'=M (L^{\frac{1}{r}}+A^{\frac{1}{r}}   +
  1)$.   Since $S_g^{-1} $ is given by $S_g^{-1}(n)=n-g(\pi_1(n))$,  we
   conclude that $S_g^{-1}$ is also  $(L', A')$ large scale Lipschitz.
  We have shown that $S_g$ is a quasiisometry.

Next we assume that  $S_g$ is a   $(L, A)$-quasiisometry, and shall prove that $g$ is large scale Lipschitz.  
   Notice that the domain of $g$ is the Euclidean space $(V_1, |\cdot|)$ (in particular,  a geodesic space). It suffices to show that there exists some constant $C\ge 0$ such that 
 $|g(u_1)-g(u_2)|\le C$ whenever $|u_1-u_2|\le 1$.

Let $u_1, u_2\in V_1$  with $|u_1-u_2|\le 1$.   
   Set  $n_1=u_1$ and $n_2=u_1*(u_2-u_1)$.   
  Then $d_h(n_1, n_2)=|u_1-u_2|\le 1$.
 Since $S_g$ is a $(L, A)$-quasiisometry,   we have $d_h(S_g(n_1), S_g(n_2))\le L+A$. 
On the other hand, we have shown that 
 $$d_h(S_g(n_1), S_g(n_2))=d_h(g(\pi_1(n_1))-g(\pi_1(n_2)),  \;  (-n_1)*n_2).$$
   Now  the generalized triangle inequality implies:
\begin{align*}
  &  d_h(0,  g(\pi_1(n_1))-g(\pi_1(n_2))) \\
& \le M\cdot \left\{d_h(0, (-n_1)*n_2))+d_h((-n_1)*n_2,\; g(\pi_1(n_1))-g(\pi_1(n_2)))\right\}\\
 & \le M\cdot (1+L+A).
\end{align*}
  Now the claim follows since 
 $\pi_1(n_1)=u_1$, $\pi_1(n_2)=u_2$ and so 
   $$d_h(0,  g(\pi_1(n_1))-g(\pi_1(n_2))) =|g(u_1)-g(u_2)|^{\frac{1}{r}}.$$

  (2)  First assume that $g$ is at finite distance from a  linear
  map  $L:  V_1\ra V_r$.   Notice that $S_L$ is an automorphism of $N$. Hence
  $S_g$ is at finite distance from the  affine map $S_L$.

  Next we assume that
$d_h(S_g, G)<\infty$, where
  $G$ is an affine map.  We shall prove that $g$ is at a finite distance from a  linear map.  
Denote by $n_0=G(0)$ and set $G_1=R_{{-n_0}}\circ G$,  where $R_{n}: \mathcal N\ra   \mathcal N$, 
 $R_{n}(x)=x*n$ is the right translation by  $n$. Then
$G_1(0)=0$. Since $d_h(G, G_1)=d_h(0,{-n_0})$, we have  $d_h(S_g,
G_1)<\infty$.
   Rewrite $G_1$  as
     $$G_1=(R_{{-n_0}}\circ L_{{n_0}})\circ
   L_{{-n_0}}\circ G,$$
  where $L_n:  \mathcal N\ra \mathcal N$ is the left translation by $n$.  
     Since $R_{{-n_0}}\circ L_{{n_0}}$ is an inner automorphism of
       $N$,  the map $G_1$ is also an affine map.
   Since $\text{Aff}(N)=N\rtimes \text{Aut}(N)$,   we see that  $G_1$ can
be written as $G_1=L_n\circ A$ for some $n\in N$ and some
automorphism $A$ of $N$. Since $G_1(0)=0$ and $A(0)=0$, we have
$n=0$. That is, $G_1$ itself is an automorphism.

Since $G_1$ is an automorphism, there exist linear maps $f_i: V_1
\ra V_i$ such that $G_1(n)=f_1(n)+\cdots  + f_r(n)$  for $n\in V_1$.
Notice that $\pi_1(n)=n$ for $n\in V_1$.
  It follows that
    $$(-S_g(n))*G_1(n)=(-n-g(n))* ( f_1(n)+\cdots +
f_r(n))$$
  for $n\in V_1$.  
     The CBH formula  yields  that the $V_r$ component
of
    $(-S_g(n))*G_1(n)$   is $P(n)-g(n)$ where $P: V_1\ra V_r$ is a  polynomial  map  of
    degree at most $r$.    It follows from    $d_h(S_g,
G_1)<\infty$ that
      $g$  is at a finite distance from   the 
polynomial map   $P$.
  Since $g$ is  a large scale Lipschitz map,
  so is   the polynomial map $P:
     V_1\ra   V_r$ . 
   But
     this is possible only when $P$ has degree one.   Hence $g$ is
     at finite distance from a linear map.

\end{proof}

From now on we shall assume that the map $g$ is large scale   Lipschitz.

We observe that   when $N$ is non-abelian,    every shear map induces the identity map on the
asymptotic cones.  To be more precise,   if the sequences of base points 
  are the constant sequences $*=(0)$,  then for any sequence 
 $\lambda_i\ra 0$  ($\lambda_i>0$) and any ultra filter $\omega$,  the induced map
   $$\tilde{S_g}: \text{Cone}_\omega(\mathcal N)\ra \text{Cone}_\omega(\mathcal N)$$
  is the identity map.      This follows from the fact that $S_g$ deviates from the identity map sublinearly:  while $g$ may grow linearly, the term $g(\pi_1(n)$ lies in  $V_r$ and the   restriction of the Carnot metric  (or $d_h$)  to $V_r$   is biLipschitz equivalent to 
   $|\cdot|^{\frac{1}{r}}$.     More concretely,  the above claim follows from the following Lemma.

\b{Le}\label{cone}  Assume $r\ge 2$.  Then 
    for any sequence  $\lambda_i\ra 0$  ($\lambda_i>0$),  the sequence of maps
  $f_i=\delta_{\lambda_i}\circ  S_g\circ \delta^{-1}_{\lambda_i}:  \mathcal N\ra \mathcal N$ converges to the identity map uniformly on compact subsets.  Here
 for any $\lambda>0$, the map $\delta_\lambda: \mathcal N\ra \mathcal N$ is the standard dilation of $\mathcal N$ given by $\delta_\lambda(v)=\lambda^j v$ for $v\in 
 V_j$. 

\end{Le}

\b{proof}
By a direct calculation we obtain:
   $f_i(n)=n+\lambda_i^r\cdot g(\frac{1}{\lambda_i} \pi_1(n))$. Hence
 \b{align*}
d_h(n, f_i(n))&=d_h(0,  \lambda_i^r\cdot g(\frac{1}{\lambda_i} \pi_1(n))\\
& =\left|\lambda_i^r\cdot g(\frac{1}{\lambda_i} \pi_1(n))\right|^{\frac{1}{r}}\\
&= \lambda_i\cdot\left| g(\frac{1}{\lambda_i} \pi_1(n)) -g(0)+g(0)\right|^{\frac{1}{r}}\\
& \le \lambda_i\cdot \left| g(\frac{1}{\lambda_i} \pi_1(n)) -g(0)\right|^{\frac{1}{r}} +\lambda_i\cdot \left|g(0)\right|^{\frac{1}{r}}\\
&\le \lambda_i\cdot \left|  L\cdot  |\frac{1}{\lambda_i} \pi_1(n)  |+A \right|^{\frac{1}{r}}   +\lambda_i\cdot \left|g(0)\right|^{\frac{1}{r}}\\
&=\left|L\cdot \lambda_i^{r-1}|\pi_1(n)|+\lambda_i^r A\right|^{\frac{1}{r}}
+\lambda_i\cdot \left|g(0)\right|^{\frac{1}{r}}\ra 0
\end{align*}
uniformly on compact subsets as $i\ra \infty$.

\end{proof}

Next we   treat the general case of   simply connected nilpotent Lie
groups.  Let $N$ be a   $r$-step  simply connected nilpotent Lie
group   with  Lie algebra $\mathcal{N}$.     The descending central
series of $\mathcal{N}$  is defined   as follows:
  $$\mathcal{N}_0=\mathcal{N}, \;  
     \mathcal{N}_{i}=[\mathcal{N},  \mathcal{N}_{i-1}],       1\le  i\le   r-1.$$
        Let
          $V_i$   be a   complementary  linear subspace  of $\mathcal{N}_i$ in
  $ \mathcal{N}_{i-1}$, $1\le i\le r-1$, and set $V_r=\mathcal{N}_{r-1}$.
    Then $\mathcal{N}$ admits the   direct sum decomposition
     $\mathcal{N}=V_1\oplus \cdots  \oplus   V_r$.
       Let $\pi_i:  \mathcal{N}  \ra  V_i$ be the projection with
       respect to this direct sum decomposition.

       Since any
two left invariant Riemannian metrics on $N$ are biLipschitz equivalent,
  we may pick  a particular such metric that is convenient for our purposes.
     Such a metric is provided by  the following   result,  see \cite{G}  or \cite{K}.

\b{theorem}\label{karidi}
     Let $N$ be a $r$-step   simply connected nilpotent Lie group with Lie algebra $\mathcal{N}$.
  Then there exist  a complementary  linear subspace  $V_i$ of $\mathcal{N}_i$ in
  $ \mathcal{N}_{i-1}$, $1\le i\le r-1$,  and an inner product on $\mathcal{N}$ such that
    all the linear subspaces $V_i$  ($1\le i\le r$, where   we set $V_r=  \mathcal{N}_{r-1}$)  are perpendicular to each other and
    such that
      the associated left invariant Riemannian   distance
          $d$  has the following property.
          There exists a constant $a>1$ such that   for each $R\ge 1$,    the ball   with radius $R$ centered at the origin
              is contained in the box  $\{n\in \mathcal{N}:
              |\pi_i(n)|\le(aR)^i  \;\text{for}\; 1\le i\le r\}$,
                and  contains the box
                  $\{n\in  \mathcal{N}:    |\pi_i(n)|\le(R/a)^i  \;\text{for}\; 1\le i\le
                  r\}$.

\end{theorem}

  For $n\in \mathcal{N}$,  set
   $$|n|_h=\sum_{i=1}^r |\pi_i(n)|^{\frac{1}{i}}.$$
     Theorem   \ref{karidi}  implies that if $d(0, n)\ge 1$, then
         $$d(0,n)/a\le |n|_h\le ar\cdot d(0,n).$$

Let  $g:  V_1\ra    V_r{}$ be any map.  
      The shear map  $S_g: \mathcal{N}\ra \mathcal{N}$  is defined  by
     $$S_g(n)=n+g(\pi_1(n)).  $$

\b{theorem}\label{main2}
  Let $S_g$ be the  shear  map   defined above.
           Then\newline
           (1) $S_g$ is a quasiisometry  from $(\mathcal{N}, d)$ to itself if and only if $g$ is large scale Lipschitz;\newline
      (2)  Let $g$ be a large scale Lipschitz map.   Then    $S_g$ is  at finite distance from an affine  map if and
      only if    $g$ is   at finite distance from a linear map.

              \end{theorem}

\b{remark}\label{rk1}
   The proof shows  that the map $S_g$ is biLipschitz if $g$ is Lipschitz. 

 \end{remark}

 \b{proof}
(1)   We only prove the \lq\lq if" part;  the \lq\lq only if" part is a slight modification of the proof  of  Theorem \ref{main1}.

The proof is similar to that of   Theorem \ref{main1} (1) (with
$d_h$ replaced by $d$), except in the estimate of
 $d(g(\pi_1(n_1)), g(\pi_1(n_2)))$.  The fact that $V_r$ is central in
 $\mathcal{N}$
 implies that the line segment in $V_r$ between $g(\pi_1(n_1))$ and
 $g(\pi_1(n_2))$  has length $|g(\pi_1(n_1))-g(\pi_1(n_2))|$ with respect
 to   the chosen left invariant Riemannian metric.  It follows that
   $$d(g(\pi_1(n_1)), g(\pi_1(n_2)))\le |g(\pi_1(n_1))-
 g(\pi_1(n_2))|.$$
  Since $g$ is a $(L, A)$ large scale Lipschitz map with respect to the Euclidean metrics, we have
$$d(g(\pi_1(n_1)), g(\pi_1(n_2)))\le  |g(\pi_1(n_1))- g(\pi_1(n_2))|\le L\cdot
  |\pi_1(n_1)- \pi_1(n_2)|+A. $$
   The \lq\lq if" part in Statement  (1)  then
  follows from the following  claim.

\noindent
  {\bf{Claim}}   There  exists  a constant $L_1\ge 0$ such that $|\pi_1(n_1)- \pi_1(n_2)|\le L_1\cdot d(n_1, n_2)$
    for all $n_1, n_2\in \mathcal{N}$.

 Now we prove the  claim. The CBH formula shows that
$\pi_1((-n_1)*n_2)=\pi_1(n_2)-\pi_1(n_1)$. Also notice that
  $d(n_1, n_2)=d(0, (-n_1)*n_2)$.  Hence it suffices to show
    $$|\pi_1((-n_1)*n_2)|\le L_1 \cdot d(0, (-n_1)*n_2).$$
First assume $d(0, (-n_1)*n_2)\le 1$.
   Then $(-n_1)*n_2$ lies in the ball with radius $1$ which is
   compact.   Since any two Riemannian metrics are biLipschitz on
   compact subsets, there is some $L_2> 0$ such that
      $|(-n_1)*n_2|\le L_2 \cdot d(0, (-n_1)*n_2)$.
        By construction the direct summands of $\mathcal{N}$
           are perpendicular  to each other,  so
   $$|\pi_1((-n_1)*n_2)|\le |(-n_1)*n_2|\le L_2 \cdot d(0,
(-n_1)*n_2). $$
     Next we assume $d(0, (-n_1)*n_2)\ge 1$.
    By  Theorem  \ref{karidi},
      $$|\pi_1((-n_1)*n_2)| \le |(-n_1)*n_2|_h\le  ar \cdot  d(0,(-n_1)*n_2).$$

(2)  The proof is  the same as that of Theorem \ref{main1} (2).

 \end{proof}

We remark that there are many
   large scale Lipschitz maps $g:      V_1\ra V_r$ that are not at finite distance from any linear maps. For the simplest examples, one can
    consider the absolute value function and power functions $|x|^\alpha$ with $0  <\alpha<1$.
  As a consequence of the constructions in this section, for any
simply connected nilpotent Lie group $N$, there are many
 quasiisometries $N\ra N$  not at finite distance from any  affine
 map.  
  Here the 
nilpotent Lie group
       $N$  is  equipped with a  left invariant Riemannian metric, and in the case of a Carnot group $N$ is  equipped  either with a Carnot metric or a left invariant Riemannian metric.

Another consequence of our construction is that there exist
quasiisometric homeomorphisms  from the quarternion Heisenberg group
 (equipped with a Carnot metric)  
to itself that are not at finite distance from any biLipschitz map.
This follows from
   Theorem \ref{main1}    and Pansu's theorem that a self biLipschitz map of the  quarternion
Heisenberg group   (equipped with a Carnot metric)  
   is a composition of left translations and  graded
automorphisms.   Notice that  the shear map   $S_g$ is a
    homeomorphism    if the map $g$ is continuous.
   One should compare this with  the    example of
  Dranishnikov-Ferry-Weinberger \cite{DFW} :
     they constructed two   
     uniformly contractible Riemannian manifolds
and a quasiisometry between them that is not at a finite distance
from any homeomorphism. The two manifolds they constructed  are
homeomorphic to the Euclidean space, but  the metrics are somehow
exotic.

As observed in the proof of   Theorem \ref{main1} (1), for  any
large scale Lipschitz map $g$, we have $S_g^{-1}=S_{-g}$.  If $g_1,
g_2: V_1\ra V_r$   are two large scale Lipschitz maps, then
$S_{g_1}\circ S_{g_2}=S_{g_1+g_2}$. It follows that the set $S$ of
all shear maps is  an abelian  group  under composition.

Let  $\text{Aff}_G(N)$ be the group of bijections of $N$ generated
by $\text{Aff}(N)$  and   $S$.
  It is easy to check that both $N$ and $\text{Aut}(N)$  normalize
  $S$, and   $\text{Aff}_G(N)=(S\cdot N)\cdot  \text{Aut}(N)$.
   We notice that the   intersections  $N\cap S$ and
   $\text{Aut}(N)\cap S$ are
    nontrivial.
      An element of
$\text{Aff}_G(N)$  will be called a \e{generalized affine map} of $N$.

\section{Other examples}\label{s2}

In this   section we construct two other examples.
  Both are biLipschitz maps of the 3-dimensional  Heisenberg group
    (equipped with a Carnot metric).
The first  biLipschitz map  sends vertical lines to 
curves  that are not vertical lines. The second   biLipschitz map   sends  vertical lines to
vertical lines and induces  the identity map  on the asymptotic
cones, but is not at finite distance from any generalized affine
map.   
     Compare 
    with Question \ref{q1}.

  Quasiisometries between  solvable Lie groups often  coarsely 
    preserve certain foliations, see \cite{EFW1}, \cite{EFW2},
        \cite{FM},
         \cite{SX}, \cite{X1}, \cite{X2}, \cite{X3}.
When studying quasiisometries between nilpotent Lie groups,
  one naturally wonders  
 whether quasiisometries coarsely  preserve cosets of
  certain special subgroups, for example,  cosets of 
  the center.      
  Our first  example shows that this is not the
  case in general.

  The  3-dimensional  Heisenberg group  $\h$ can be identified with
$\C\times \R$   with group operation given by
  $$(z_1, t_1)*(z_2, t_2)=(z_1+z_2, t_1+t_2+2 \,{\text{Im}}(z_1
  \bar{z}_2)).$$
    Let  $\pi:  \h=\C\times \R\ra \C$ be the natural projection.
  Then $\pi$ is $1$-Lipschitz with respect to the  usual Carnot metric on $\h$ and the Euclidean metric on $\C$. 
  The Heisenberg inversion $j: \h\backslash\{(0,0)\}\ra \h\backslash\{(0, 0)\}$ is the map given by (see  \cite{KR} page 315):
   $$j(z,t)=\left(\frac{z}{it-|z|^2},\;  \frac{-t}{t^2+|z|^4}\right).$$
      Recall that (the one point compactification of)  $\h$ can be
      identified
        with the ideal boundary of the complex hyperbolic plane:  $\partial\h^2_{\C}=\h\cup \{\infty\}$.
            The Heisenberg inversion corresponds to the antipodal
            map in the ball model for the complex hyperbolic plane.
  Let $\lambda>1$ be a constant.
  Let $f_\lambda: \h\ra \h $ be the map given by
   $$f_\lambda(x+iy, t)=(\lambda x+i\frac{1}{\lambda}y, t). $$
     The sought map  $F_\lambda: \h\ra \h$ is given by 
$F_\lambda(0,0)=(0,0)$  and  
$F_\lambda(p)= j  \circ f_\lambda\circ   j(p)$  for $p\in \h\backslash\{(0,0)\}$.

\b{Le}\label{bilip}
  The map $F_\lambda$ is   biLipschitz. 
\end{Le}

\b{proof}
Recall that (see \cite{SX},  Proposition 5.4 for more details)  the relation between the metrics on   $\partial  \h^2_\C$
and on $\h$ is similar to the relation between the spherical metric on the sphere and the
Euclidean metric on the Euclidean space. 
More precisely, if  $d_0$ denotes the Carnot metric
on  $\partial  \h^2_\C=\h\cup  \{\infty\}$, then the Carnot metric  $d$ on 
  $\h$   is biLipschitz equivalent to
\b{equation}\label{e1}
d_1(x,  y) = \frac{d_0(x, y)}{
d_0(\infty,  x) d_0(\infty,  y)}   \;\;\text{for}\;\; x,y\in \h.
\end{equation}
It follows that every biLipschitz map of
 $(\h\cup  \{\infty\}, d_0)$ that fixes $\infty$  
restricts to a biLipschitz map of  $(\h, d)$.
  Furthermore,    $d_0$ is biLipschitz equivalent to the function 
 $d_2:  \partial  \h^2_\C\times \partial  \h^2_\C \ra [0, \infty)$ defined by:
 $d_2(\infty, \infty)=0$, 
\b{align}
d_2(x, \infty)=d_2(\infty, x)=\frac{1}{1+d(x,0)}\;\; \text{for}\;\; x\in \h \label{e2} \\
d_2(x,y)=\frac{d(x,y)}{[1+d(x,0)][1+ d(y,0)]} \;\;\text{for}\;\; x,y\in \h\label{e3}.
\end{align}
  It follows that 
if $f: (\h, d) \ra (\h, d)$ is biLipschitz  and $f(0)=0$, then $f$ extends (by setting 
$f(\infty)=\infty$)  to  a biLipschitz map of  $(\h\cup  \{\infty\}, d_0)$.

Now since 
 $f_\lambda$  is a biLipschitz map of  $\h$  and $f_\lambda(0)=0$, it extends (with
 $f_\lambda(\infty) = \infty$) to a biLipschitz map of the ideal boundary 
 $(\partial  \h^2_\C, d_0)$.    Being the antipodal map
in the ball model,   $j$ is clearly a biLipschitz map of  $\partial  \h^2_\C$. 
  Hence  $F_\lambda$  is a biLipschitz map of  $(\partial  \h^2_\C, d_0)$. 
    Since $F_\lambda(\infty) = \infty$, 
  we see that as a map of  $(\h, d)$,    $F_\lambda$  is also biLipschitz.

\end{proof}

We next show that $F_\lambda$ is not at finite distance from any map that sends vertical lines to vertical lines.  For any $z\in \C$,   denote by $q_z$ the vertical line above $z$:
  $q_z=\pi^{-1}(z)$.  

\b{Le}\label{property.of.F}
(1)   For any $z\in \C$,  $F_\lambda(q_z)$ lies at a finite distance from a vertical line.\newline
 (2)  There is no uniform bound  for the distance  between $F_\lambda(q_z)$ (as $z\in \C$ varies) and vertical lines; in other words,
   $$\sup_{z\in \C} \inf_{w\in \C}d_H(F_\lambda(q_z), q_w)=\infty,$$
 where $d_H$ denotes the  Hausdorff distance.\newline
  (3)  $F_\lambda$ is not at finite distance from any map that sends vertical lines to vertical lines.\newline
  (4)  $\delta_a\circ  F_\lambda= F_\lambda\circ \delta_a$ for all $a>0$; 
 in particular,  the map induced by $F_\lambda$  on asymptotic cones is the map $F_\lambda$ itself. 

\end{Le}

\b{proof}
  (1)   We shall show  $d_H(F_\lambda(q_z), q_{f_\lambda^{-1}(z)})<\infty$.
Fix  $z=x+iy\in \C$.  For any $t\in \R$,     write 
 $j(z, t)=(x_1, y_1, t_1)$,  $f_\lambda\circ j(z,t)=(x_2, y_2, z_2)$ and
 $F_\lambda(z,t)=(x_3, y_3, t_3)$.  
  Then 
\b{equation}\label{e5}
x_1=\frac{yt-x|z|^2}{t^2+|z|^4},  \;\;\; y_1=\frac{-(xt+y|z|^2)}{t^2+|z|^4},  \;\;\;
t_1=\frac{-t}{t^2+|z|^4}.  \;\;\;  
\end{equation}
We shall view $x_i, y_i, t_i$  as functions of $t$.  Using (\ref{e5}) and the definition of 
  $f_\lambda$,  we can write $x_3$, $y_3$, and $t_3$ as rational functions   of $t$.
 The complete formulas  are  a little too complicated to write down here, so we shall only 
 record  the leading 
terms in the  numerators and the denominators. 
  For $x_3$, they are  $\frac{x}{\lambda}t^6$ and $t^6$.  
 So  $x_3\ra \frac{x}{\lambda}$ at $|t|\ra \infty$.
  Similarly,  $y_3\sim\frac{\lambda y t^6}{t^6}\ra \lambda y$
  and $t_3\sim \frac{t^7}{t^6}= t$ as  $|t|\ra \infty$. 
    Notice that $f_\lambda^{-1}(z)=\frac{x}{\lambda}+i\lambda y$. 
     So as $|t|\ra \infty$,
 the curve  $F_\lambda(q_z) $ is asymptotic to $q_{f_\lambda^{-1}(z)}$. Hence
  $d_H(F_\lambda(q_z), q_{f_\lambda^{-1}(z)})<\infty$.
  
 (2) Notice that for any $w\in \C$ and any $a>0$,   the $a$-neighborhood of $q_w$ is simply  $\pi^{-1}(B(w, a))$,   where $B(w,a)\subset \C$ is the ball  with center $w$ and radius $a$.    Hence it suffices to show that there is no uniform bound for 
 the diameter of $\pi(F_\lambda(q_z))\subset \C$   as   $z\in \C$  varies.
    Let $z_n=n+in$.  A direct calculation using (\ref{e5}) and the definition of 
 $f_\lambda$  shows  that   $F_\lambda(z_n, 0)$ has $y$-coordinate
  $\frac{2\lambda}{\lambda^4+1}n$.  On the other hand, in the proof of (1) we saw that 
  the $y$-coordinates of  the points on $F_\lambda(q_{z_n})$  converge to $\lambda n$. 
 Since $\lambda>1$,  we have  $\lambda n-\frac{2\lambda}{\lambda^4+1}n\ra \infty$ as $n\ra \infty$.      Hence there is no 
uniform bound for 
 the diameter of $\pi(F_\lambda(q_z))$.

(3) follows  easily  from  (2).

 (4)   Notice that $\delta_a\circ j\circ \delta_a=j$  and 
  $\delta_a\circ  f_\lambda= f_\lambda  \circ \delta_a$ for any $a>0$.
  It follows that 
 \b{align*}
\delta_a\circ  F_\lambda&=\delta_a\circ 
(j\circ f_\lambda\circ j)\\
 &=(\delta_a\circ j\circ \delta_a)\circ (\delta_{\frac{1}{a}}\circ f_\lambda\circ \delta_a)\circ 
(\delta_{\frac{1}{a}}\circ j\circ \delta_{\frac{1}{a}})\circ  \delta_a\\
&=(j\circ f_\lambda\circ j)\circ  \delta_a\\
&=F_\lambda\circ  \delta_a.
\end{align*}
In particular,  $\delta_a\circ  F_\lambda\circ \delta_{\frac{1}{a}}=F_\lambda  $ converges to $F_\lambda$ as 
$a\ra 0$.   The claim follows. 

\end{proof}

It is an open question whether there exist biLipschitz maps of the Heisenberg group that map
  vertical lines to curves at infinite Hausdorff distance from vertical lines.

\b{remark}\label{r5}
  In the construction of $F_\lambda$, we can replace $f_\lambda$ with any 
 biLipschitz map $f: \h\ra \h$ that fixes $0$.  For a general 
 biLipschitz map $h: \h\ra \h$,  if $p:=h(0)$, then the map $f:=L_{-p}\circ h:\h\ra \h$ is also 
 biLipschitz   and satisfies  $f(0)=0$.    See Theorem \ref{tyson} below for more examples of biLipschitz maps of  $\h$. 

\end{remark}

 We next construct our second  biLipschitz map of $\h$.  
Let $f: \C\ra \C$ and $F: \h\ra \h$ be two maps.  We say $F$ is a  
   lift of $f$ if 
 $f\circ \pi=\pi\circ F$.  
 Below we shall identify $\C$ with $\R^2$.  
 We  shall use 
 the following result:

\begin{theorem}\label{tyson}
\e{(Balogh, Hoefer-Isenegger,  Tyson \cite{BHT},  Theorem 1.6)}\newline
  (1)  Let $f:  \R^2\ra \R^2$ be Lipschitz with ${\text{det}} (Df)= \lambda$ almost everywhere. Then  there exists a Lipschitz lift
    $F:  \h\ra  \h$.  
       If  $\tilde F$ is another  Lipschitz lift of $f$, then $\tilde{F}(z, t)=F(z,t+\tau)$ for some $\tau\in \R$.  Here $\text{det}$  denotes the determinant and $Df$  denotes  the differential of $f$. 
  \newline
        (2)  Conversely,
             if $f:\R^2\ra \R^2$  is Lipschitz with Lipschitz lift, then there exists $\lambda\in \R$ so that
             ${\text{det}} (Df)= \lambda$ almost everywhere.
             \end{theorem}

 Let  $f: \R^2\ra \R^2$    be  Lipschitz  
with ${\text{det}} (Df)= \lambda$ almost everywhere.       Write $f$  
as $f(x,y)=(f_1(x,y), f_2(x,y))$,    and   denote by 
   $Df^*$   the transpose of $Df$. 
              It is shown in \cite{BHT} that  a Lipschitz lift  $F$   is given by
               $$F(z,t)=(f(z), \lambda t+ h_0(z))$$
                 where  $h_0: \C=\R^2\ra \R$ is determined by 
                 $$  \begin{pmatrix} \frac{\partial h_0}{\partial x}(x,y)\\
\frac{\partial h_0}{\partial y}(x,y)\end{pmatrix}
=\begin{pmatrix}-2\lambda y\\2\lambda x\end{pmatrix}-2\cdot Df^*\cdot \begin{pmatrix}-f_2(x,y)\\ f_1(x,y) \end{pmatrix}\;\; a.e. $$

Let $g: \R\ra \R$ be the map $g(y)=y$ for $|y|\le 1$ and
$g(y)=y^{\frac{1}{3}}$ for $|y|\ge 1$.
  Let $f: \C\ra \C$ be given by $f(x+iy)=(x+g(y))+iy$.
  Then  $f$ is biLipschitz   with ${\text{det}}(Df)\equiv 1$. 
Hence by
  Theorem \ref{tyson},
$f$ has a biLipschitz lift $F: \h\ra \h$. The map $F$ is given by
    $$F(x+iy, t)=(x+g(y)+iy, t+h_0(y)),$$
         where
  $h_0(y)=-1$ for $|y|\le 1$ and $h_0(y)=-y^{\frac{4}{3}}$ for
  $|y|\ge 1$.

We list some of the properties of the map $F$ in the following  Lemma.

\b{Le}\label{e3}
(1)  $F$ maps vertical lines to vertical lines;\newline
(2) $F$ induces the identity map on asymptotic cones;\newline
 (3)   $F$ is not at finite distance from any
generalized affine map.

\end{Le}

\b{proof}
(1) is clear from the definition of $F$.

(2)   For any $a>0$, the map $\delta_{a}\circ  F\circ \delta^{-1}_{a}$  is given by:
  $$\delta_{a}\circ  F\circ \delta^{-1}_{a}(x+iy,\; t)=\left(x+a\cdot g\left(\frac{y}{a}\right)+iy, \;t+a^2\cdot h_0\left(\frac{y}{a}\right)\right).$$
     Now it is easy to see that 
    for any sequence  $a_i\ra 0$  ($a_i>0$),  the sequence of maps
  $F_i:=\delta_{a_i}\circ  F\circ \delta^{-1}_{a_i}:  \h\ra \h$ converges to the identity map uniformly on compact subsets.   Hence $F$ induces the identity map 
 on asymptotic cones.  

 (3)  Notice that a generalized affine map  of $\h$ is a lift of an affine map of $\R^2$. 
 Also recall that $\pi: \h\ra \R^2$ is $1$-Lipschitz.  Suppose $F$ is 
at finite distance from   some  
generalized affine map. Then it follows that $f$ is at finite distance from an affine map
  of $\R^2$,  which is not true.
\end{proof}

The theorem of Balogh, Hoefer-Isenegger    and Tyson shows that
Lipschitz maps of the Heisenberg group preserving the vertical
foliation are very restricted.  
    Is there a similar result for other Carnot groups?    Another  interesting  question is whether there is an
 analogous result for quasiisometries of the Heisenberg group
 that coarsely preserve the vertical foliation.

\section{Open questions}\label{s3}

We end this paper by raising a few questions about the structure of quasiisometries between nilpotent Lie groups.

For the quarternion   Heisenberg group,
   a natural question  is the following (compare with the second example in Section \ref{s2}):

    \b{question}\label{q1}
 Is every  quasiisometry of the quarternion  Heisenberg group
    at finite distance from a generalized affine map?
\end{question}

Let $N$ be a simply connected nilpotent Lie group,  and $N_1=N$,   $N_2$, $\cdots$,   $N_r$ be the subgroups in the lower central series of $N$.    

 \b{question}\label{q2}
  Let $N$ be   a    simply   connected
 nilpotent Lie group equipped with a left invariant Riemannian metric or a 
Carnot group  equipped with  a  Carnot metric. 
     Suppose the asymptotic cone of $N$ does not have a Heisenberg group factor.  
  Is it true that every quasiisometry  $N\ra N$      coarsely    permutes the left cosets of  
  $N_i$ for each $i$?

\end{question}

 Notice that the answer to the above question is negative in the presence  of a 
Heisenberg group factor:  for any $N$, the map 
 $$id\times F_\lambda: N\times \h\ra N\times \h$$
     is a counterexample,  where $F_\lambda$ is the  map in the  first  example of  Section \ref{s2}. 

On the level of asymptotic cones,  Question \ref{q2} becomes:

\b{question}\label{q3}
 Let  $N$ be a Carnot group equipped with a Carnot metric.  Assume $N$  is not isomorphic to the direct product  of  an  Heisenberg group  and another  Carnot group.  Is it true that  every biLipschitz map $F: N\ra N$   permutes the left cosets of  
  $N_i$ for each $i$?

\end{question}

 \addcontentsline{toc}{subsection}{References}

\noindent Xiangdong Xie: Dept. of Mathematics and Statistics,  Bowling Green State University,   Bowling Green, OH 43403, U.S.A.\hskip .4cm E-mail:
  xiex@bgsu.edu

\end{document}